\date{}

\date{}

\title{A factor matching of optimal tail between Poisson processes}
\author{\'Ad\'am  Tim\'ar }

\documentclass[11pt,naturalnames]{article}
\renewcommand\footnotemark{}
\usepackage[ruled,vlined]{algorithm2e}
\usepackage{amsmath}
\usepackage{amsthm}
\usepackage{amsfonts}
\usepackage{graphicx}
\newif\ifhyper\IfFileExists{hyperref.sty}{\hypertrue}{\hyperfalse}
\ifhyper\usepackage{hyperref}\fi
\usepackage{enumitem}
\usepackage{subfig}
\usepackage{caption}
\usepackage{keyval}
\usepackage[margin=1in]{geometry}
\usepackage{setspace}

\usepackage[displaymath, mathlines,pagewise]{lineno}

\theoremstyle{definition}

\newtheorem{theorem}{Theorem}

\newtheorem{remark}[theorem]{Remark}

\newtheorem{definition}{Definition}

\def \proof {{ \medbreak \noindent {\bf Proof.} }}
\def\proofof#1{{ \medbreak \noindent {\bf Proof of #1.} }}
\def\proofcont#1{{ \medbreak \noindent {\bf Proof of #1, continued.} }}
\def\supp{{\rm supp}}
\def\max{{\rm max}}
\def\min{{\rm min}}
\def\dist{{\rm dist}}
\def\Aut{{\rm Aut}}
\def\id{{\rm id}}
\def\Stab{{\rm Stab}}

\begin{document}
\maketitle
\let\thefootnote\relax\footnotetext{\footnotesize{Partially supported by Icelandic Research Fund Grant 185233-051 and the ERC Consolidator Grant 772466 ``NOISE''.}}

\bigskip

\def\eref#1{(\ref{#1})}
\newcommand{\Prob} {{\bf P}}
\newcommand{\C}{\mathcal{C}}
\newcommand{\LL}{\mathcal{L}}
\newcommand{\Z}{\mathbb{Z}}
\newcommand{\N}{\mathbb{N}}
\newcommand{\HH}{\mathbb{H}}
\newcommand{\Rr}{\mathbb{R}^3}
\newcommand{\h}{\mathcal{H}}
\def\diam{\mathrm{diam}}
\def\length{\mathrm{length}}
\def\ev#1{\mathcal{#1}}
\def\Isom{{\rm Isom}}
\def\Re{{\rm Re}}
\def \eps {\epsilon}
\def \P {{\Bbb P}}
\def \E {{\Bbb E}}
\def \proof {{ \medbreak \noindent {\bf Proof.} }}
\def\proofof#1{{ \medbreak \noindent {\bf Proof of #1.} }}
\def\proofcont#1{{ \medbreak \noindent {\bf Proof of #1, continued.} }}
\def\supp{{\rm supp}}
\def\max{{\rm max}}
\def\min{{\rm min}}
\def\dist{{\rm dist}}
\def\Aut{{\rm Aut}}
\def\id{{\rm id}}
\def\Stab{{\rm Stab}}

\newcommand{\lra}{\leftrightarrow}
\newcommand{\xlra}{\xleftrightarrow}
\newcommand{\xnlra}{\xnleftrightarrow}
\newcommand{\pc}{{p_c}}
\newcommand{\pt}{{p_T}}
\newcommand{\ptk}{{\hat{p}_T}}
\newcommand{\pl}{{\tilde{p}_c}}
\newcommand{\pe}{{\hat{p}_c}}
\newcommand{\pr}{\mathrm{\mathbb{P}}}
\newcommand{\pp}{\mu}
\newcommand{\ex}{\mathrm{\mathbb{E}}}
\newcommand{\ee}{\mathrm{\overline{\mathbb{E}}}}

\newcommand{\om}{{\omega}}
\newcommand{\ebd}{\partial_E}
\newcommand{\ivbd}{\partial_V^\mathrm{in}}
\newcommand{\ovbd}{\partial_V^\mathrm{out}}
\newcommand{\q}{q}
\newcommand{\TT}{\mathfrak{T}}
\newcommand{\T}{\mathcal{T}}
\newcommand{\RR}{\mathcal{R}}

\newcommand{\CC}{\Pi}
\newcommand{\BB}{\Pi}

\newcommand{\A}{\mathcal{A}}
\newcommand{\cc}{\mathbf{c}}
\newcommand{\pa}{{PA}}
\newcommand{\degi}{\deg^{in}}
\newcommand{\dego}{\deg^{out}}
\def\Pn{{\bf P}_n}

\newcommand{\R}{\mathbb R}
\newcommand{\F}{F}
\newcommand{\FF}{\mathfrak{F}}
\newcommand{\Can}{\rm Can}
\newcommand{\Vol}{\mathrm Vol}

\def\Gstar{{{\cal G}_{*}}}
\def\Gstarstar{{{\cal G}_{**}}}
\def\Gstarplus{{{\cal G}_{*}^\frown}}
\def\Rel{{\cal R}}
\def\Comp{{\rm Comp}}
\def\calG{{\cal G}}
\def\omps{{\omega_\delta^\eps}}

\begin{abstract}
Consider two independent Poisson point processes of unit intensity in the Euclidean space of dimension $d$ at least 3. We construct a perfect matching between the two point sets that is a factor (i.e., an equivariant measurable function of the point configurations), and with 
the property that the distance between a configuration point and its pair has a tail distribution that decays as fast as possible, namely, as $b\exp (-cr^d)$ with suitable constants $b,c>0$.
Our proof relies on two earlier results: an allocation rule of similar tail for a Poisson point process, and a recent theorem that enables one to obtain perfect matchings from fractional perfect matchings in our setup.
\end{abstract}

\bigskip

Let $\omega_1$ and $\omega_2$ be two independent random sets of points in the Euclidean space $\R^d$, given by Poisson point processes of unit intensity. We are interested in measurable perfect matchings (bijections) between $\omega_1$ and $\omega_2$. We want the perfect matching to commute with translations of $\R^d$, in other words, to be invariant under translations. One may allow the use of extra randomness in a version of this problem, as in \cite{HPPS}, but in our setting the matching has to be a deterministic, equivariant function of the point configurations. Such a matching is called a {\it factor matching}. 

Our goal is to make matched points ``as close as possible'', that is, we want $\P\bigl(|m(0)|>r\bigl| 0\in\omega_1\bigr)$ decay as fast as possible, where $m(0)$ is the point matched to $0$. Note that conditioning on the event (of 0 probabity) that $0\in\omega_1$ is a standard and natural construction, called the Palm version, and in case of a Poisson point process it is the same as adding a configuration point to 0; see e.g. \cite{Th} for general background.

A trivial lower bound for the possible tail is given by the closest point of $\omega_2$ to 0, whose tail distribution is of order $\exp (-cr^d)$.
Our main result is that a decay of this magnitude is attainable (up to the value of the constant).

\begin{theorem}\label{main}
Let $\omega_1$ and $\omega_2$ be independent Poisson point processes of intensity 1 in the Euclidean space $\R^d$ of dimension $d\geq 3$. There exists a factor perfect matching $m$ between $\omega_1$ and $\omega_2$ with the property that
\begin{equation}\label{fontos}
\P\bigl(|m(0)|>r\bigl|0\in\omega_1\bigr)<b\exp (-cr^d)
\end{equation}
with some $b,c>0$, where
$m(0)$ is the matched pair of 0. 
\end{theorem}
 
There are several reasons why the requirement of no extra randomness in a factor matching is of importance. It is certainly natural: configuration points have to find their pairs using only local information (up to arbitrarily small error) and no ``central planning''. Analogous problems for graphs have shown striking difference between parameter values optimized over all invariant or over all equivariant objects in cases, the most famous example possibly being the independent set of transitive graphs, where the maximal factor (of iid) independent set may have lower density than the maximal invariant independent set \cite{B}. Similar difference can be observed for 
certain functions from point processes. For example, the fastest possible decay of a perfect matching over a single Poisson point process in dimension 1 is very different depending on weather we allow extra randomness or not; see \cite{HPPS} for more details. 
The requirement of no extra randomness usually makes the problem more difficult, as illustrated by optimal allocations for a Poisson point process (which are closely related to matchings, as also seen below). A (randomized) allocation of optimal decay with the use of extra randomness was found in \cite{HP}, while the construction of factor allocations needed significant further efforts, including the analysis of gravitational allocation \cite{CPPR}, and later a factor allocation of optimal tail \cite{MT}. Similarly, for matchings between two Poisson point processes, an optimal one using extra randomness was found in the seminal work \cite{HPPS} by Holroyd, Pemantle, Peres and Schramm, while the same paper only contains a factor matching of worse than polynomial decay.
This bound was improved in a never revised preprint \cite{T} by the present author, where a matching scheme of decay $<\exp (-cr^{d-2-\eps})$ was constructed. (That paper was dealing with the essentially equivalent problem of flipping a fair coin for each vertex of $\Z^d$, $d\geq 3$, and trying to find a factor matching of optimal tail between vertices with a head and vertices with a tail.)
The case of $d=1,2$ shows a different behavior, with $c'r^{-d/2}$ being the optimal decay, as shown in \cite{M} ($d=1$) and the simpler part of \cite{T} ($d=2$).
Various related versions of the problem have been studied, such as optimal tail behavior of matchings allowing extra randomness or factor matchings on a unique Poisson point process \cite{HPPS}, minimal planar matchings \cite{H}, tail behavior of the so-called stable matching \cite{HHP} or
multicolor matchings \cite{AAH}... Let us mention that matching two Poisson point processes is significantly more difficult than the problem of finding a perfect matching factor for a single Poisson point process, essentially because of the discrepancy in the number of configuration points from the two processes within any fixed box.


\begin{definition}
Consider a translation invariant point process $\omega$ of intensity 1 in the Euclidean space $\R^d$. A {\it (factor) allocation} is a measurable map $\alpha_\omega$ from $\R^d$ to $\omega$ such that for every $x\in\omega$, $\alpha^{-1}_\omega(x)$ has unit area, and
$\alpha_\omega$ is equivariant (with regard to translations) and measurable. We will assume for simplicity that $\alpha_\omega (x)=x$ for every $x\in\omega$.  The sets $\alpha^{-1}_\omega(x)$ ($x\in\omega$) will be called {\it cells}. Whenever an allocation $\alpha_\omega$ is given, we will denote the cell $\alpha^{-1}_\omega \alpha_\omega(p)$ of an arbitrary $p\in\R^d$ by $\A (p)$.
\end{definition}

The proof will rely on two earlier results. We denote the Euclidean ball of radius $r$ around $0$ by $B(r)$.

\begin{theorem}\label{MarkoT}[Mark\'o-T., \cite{MT}]
Let $\omega$ be a Poisson point process of intensity 1 in $\R^d$, $d\geq 3$. Then there is an allocation rule with the property that 
\begin{equation}
\P\bigl(\diam (\A(0))>r\bigl|0\in\omega\bigr)<b_0\exp (-c_0 r^d),
\end{equation}
with some $b_0,c_0>0$. The allocation has the further property that any bounded subset of $\R^d$ is intersected by only finitely many allocation cells.
\end{theorem}
The last assertion of the theorem is not stated explicitly in \cite{MT} but is clear from the construction: to be in a cell of large diameter, one has to be contained in some box with an ``extra large'' number of configuration points, and the probability that there are infinitely many such boxes intersecting a given bounded set is 0.

A {\it perfect fractional matching} on a graph $G=(V(G),E(G))$ is a map $f:E(G)\to [0,\infty)$ such that $\sum_{v\in e, e\in E(G)} f(e)=1$ for every $v\in V(G)$.

\begin{theorem}\label{BKS}[Bowen-Kun-Sabok, \cite{BKS}]
Suppose that a hyperfinite one-ended bipartite graphing $G$ has some measurable perfect fractional matching which is almost everywhere positive. Then $G$ has a measurable perfect matching almost everywhere.
\end{theorem}

\begin{proofof}{Theorem \ref{main}}
First we construct a hyperfinite one-ended bipartite graphing $G$ on $\omega_1$ and $\omega_2$ as classes of bipartition, and furthermore show that there exists an almost everywhere positive measurable perfect fractional matching on $G$. 
Let $\alpha_1$ and $\alpha_2$ be the allocation functions for $\omega_1$ and $\omega_2$ respectively, as given by Theorem \ref{MarkoT}, and let $\A_1=\alpha_1^{-1}\alpha_1$ (respectively, $\A_2=\alpha_2^{-1}\alpha_2$) be the map assigning its cell to every point of $\R^d$. Let $x\in\omega_1$ and $y\in\omega_2$ be adjacent in $G$ if $\A_1(x)\cap\A_2(y)$ has positive Lebesgue measure. (This graph applied for the construction of a matching from an allocation, already appears in \cite{AKT}, in a finite setting.) Define the perfect fractional matching by assigning the measure of this intersection to the edge $\{x,y\}\in E(G)$. This will indeed be a fractional perfect matching, by definition of an allocation. The graphing $G$ that we have just defined is bipartite by definition, every vertex has finite degrees by Theorem \ref{MT}, and $G$ has one end because a finite $U\subset V(G)$ can separate $x,y\in V(G)$ only if $\cup_{u\in U} \A_{i_u}(u)$ separates the cells of $x$ and $y$ in $\R^d$
(here $i_u$ is 1 if $u\in\omega_1$, otherwise it is 2), but the removal of such a bounded set from $\R^d$ can only leave one unbounded connected component. 

So it only remains to prove that $G$ is hyperfinite, that is, to find for any $\eps>0$ a subset of the configuration points as a factor (measurable equivariant function), such that the density of the chosen points in $V(G)$ is less than $\eps$
and their removal breaks $G$ into all finite components.
Given $\eps>0$, choose $r$ such that $\P\bigl(\max \{|x|: \{0,x\}\in E(G)
\}>r\bigl| 0\in\omega_1\bigr)<\eps/2$. Let $U_1$ be the set of vertices $u$ in $V(G)$ such that $\max \{|u-x|: \{u,x\}\in E(G)
\}>r$. Then $\P(0\in U_1| 0\in\omega_1)<\eps/2$. Now choose $N$ large enough so that $2^d r/N<\eps/2$.
We will take a partition ${\mathcal P}$ of $\R^d$ into convex polyhedra as a factor of $\omega_1\cup\omega_2$, in such a way that every polyhedron in ${\mathcal P}$ contains a ball of radius $N$. This condition implies that $\P(0\in\cup_{P\in{\mathcal P}} \partial_r P)<2^d r/N<\eps/2$, with $\partial_r P$ denoting the $r$-neighborhood of the boundary of $P$. 
Hence, denoting by $U_2$ the set of vertices in $V(G)\setminus U_1$ whose cell intersects the boundary of some $P\in{\mathcal P}$, we have $\P(0\in U_2|0\in\omega_1)<\eps/2$. The removal of $U_1\cup U_2$ splits $G$ into only finite components, because points of $\omega_1\cup\omega_2$ that are in different elements of ${\mathcal P}$ cannot be in the same component of $G\setminus(U_1\cup U_2)$. Furthermore, $\P(0\in U_1\cup U_2|0\in\omega_1)<\eps$, proving hyperfiniteness. So it only remains to construct ${\mathcal P}$. To do so, it is enough to find a factor subset $X$ of $\omega_1\cup\omega_2$ consisting of points whose pairwise distance in $\R^d$ is at least $2N$, and then take the Voronoi tessellation for $X$ to be ${\mathcal P}$. One way to do this is to define a graph $H$ on $\omega_1\cup\omega_2$ where two points are adjacent if their distance is at most $2N$, then delete points whose degree is at least $D$ for a suitable $D$ so that we do not delete everything, and finally 
use a measurable Brooks theorem to find a proper coloring of the resulting induced graph by finitely many colors (see 
\cite{KST} for the measurable Brooks theorem, which is applied a few times in \cite{BKS} in a similar way as here). Choose one of the color classes to be $X$.

Now consider the graphing $G$  and let $m$ be a perfect matching that bijectively assigns to every $x\in\omega_1$ a point $m(x)\in \omega_2$, as provided by Theorem \ref{BKS}. We only have to show \eqref{fontos}. Conditional on $0\in\omega_1$, the random marked graph $(G,0;\A)$ is unimodular, where by $\A$ we denote the assignment of the number $\diam(\A_i(x))$ (with $i=1$ if $x\in\omega_1$ and  $i=2$ otherwise) to each vertex $x$. Define the following simple mass transport: let $y\in\omega_2$ send mass 1 to $x\in\omega_1$ if $m(x)=y$ and $\diam (\A_2(y))>r$. By the Mass Transport Principle, 
$$\P\bigl(\diam (\A_2(m(x)))>r
\bigl|0\in\omega_1\bigr)=\E(\text{``mass received''})
$$
$$=\E(\text{``mass sent out''})
=\P\bigl(\diam (\A_2(0))>r
\bigl|0\in\omega_2\bigr).
$$
We conclude that
$$
\P\bigl(|m(0)|>2r\bigl|0\in\omega_1\bigr)<\P\bigl(\diam (\A_1(0))>r
\bigl|0\in\omega_1\bigr)+\P\bigl(\diam (\A_2(m(0)))>r
\bigl|0\in\omega_1\bigr)<2b_0\exp (-c_0 r^d)
$$
by Theorem \ref{MarkoT}.
\end{proofof}
\qed

\medskip

\begin{remark}
In Theorem \ref{main} we construct a matching that is equivariant with respect to {\it translations}.  However, the same construction would work for all isometries, if the allocation considered from Theorem \ref{MarkoT} were isometry-equivariant. While the proof in \cite{MT} is elaborated only for translations, a bit of extra work would make it isometry-equivariant (see Remark 3.9 in \cite{MT}), giving rise to an isometry-equivariant matching in our Theorem \ref{main}.
\end{remark}

\bigskip
\medskip

\noindent
{\bf Acknowledgments:} 
I am greatly indebted to G\'abor Kun for telling me about the results in \cite{BKS} and for further discussions. I thank G\'abor Pete for his comments and suggestions.

\ \\

{\bf \'Ad\'am Tim\'ar}\\
Division of Mathematics, The Science Institute, University of Iceland\\
Dunhaga 3 IS-107 Reykjavik, Iceland\\
and\\
Alfr\'ed R\'enyi Institute of Mathematics\\
Re\'altanoda u. 13-15, Budapest 1053 Hungary\\
\texttt{madaramit[at]gmail.com}\\

\end{document}